\documentclass[12pt,a4paper]{article}

\usepackage{amssymb, amsmath, amsthm}

\newcommand{\dd}{\mathrm{d}}

\theoremstyle{plain}
\newtheorem*{theorem}{Theorem}

\theoremstyle{remark}
\newtheorem{remark}{Remark}

\title{On a conjecture of Seneta}

\author{P\'eter Kevei\thanks{
Bolyai Institute, 6720, Aradi v\'ertan\'uk tere 1. Szeged, Hungary;
e-mail address: \texttt{kevei@math.u-szeged.hu}} \\
University of Szeged
}

\date{}

\begin{document}

\maketitle

\begin{abstract}
In this short note we prove that
$h_\beta(x) = \beta \int_0^x y^{\beta-1} \overline F(y) \dd y$
is regularly varying with index $\rho \in [0,\beta)$ if and only if 
$V_\beta (x) = \int_{[0,x]} y^\beta \dd F(y)$ is regularly varying
with the same index. This implies an extended version of a recent
conjecture by Seneta \cite{Seneta19}.

\end{abstract}

\section{Introduction and results}

Let $F$ be the distribution function of a nonnegative random variable,
and put $\overline F(x) = 1 -F(x)$.  
For $\beta > 0$ introduce the truncated $\beta$-moments as
\begin{equation*} 
\begin{split}
& h_\beta(x) = \beta \int_0^x y^{\beta-1} \overline F(y) \dd y \\
& V_\beta (x) = \int_{[0,x]} y^\beta \dd F(y).
\end{split}
\end{equation*}
Integration by parts gives 
\begin{equation} \label{eq:int-parts}
\begin{split}
V_\beta(x)
& = 
\beta \int_0^x y^{\beta-1}  \overline F(y) \dd y - 
x^{\beta}  \overline F(x) \\
& = h_\beta (x) - x^{\beta} \overline F(x),
\end{split}
\end{equation}
the basic relation between $h_\beta$ and
$V_\beta$. 
Note that if $h_\beta$ or $V_\beta$ is regularly varying then
its index $\rho$ is at most $\beta$. 

If $\int_{[0,\infty)} y^\beta \dd F(y) < \infty$ then, as $x \to \infty$,
$x^\beta \overline F(x) \to 0$ and
$h_\beta(x) \to V_\beta(\infty)$. Therefore, we always
assume
\begin{equation} \label{eq:beta-ass}
\int_{[0,\infty)} y^\beta \dd F(y) = \infty.
\end{equation}

Seneta \cite[Theorem 3]{Seneta19} proved among other equivalences
that if $V_\beta$ is slowly
varying then $h_\beta$ is slowly varying and conjectured 
\cite[p.1653]{Seneta19} that the converse also holds.
In this short note we prove an extension of Seneta's conjecture.
The following result is a generalization of Theorem 3 
by Seneta \cite{Seneta19}, where he considers 
the case $\rho  = 0$.

In what follows $\mathcal{RV}_\rho$ stands for the class of regularly
varying functions with index $\rho$, and any nonspecified limit 
relations are meant as $x \to \infty$.

Before the main result we need some definitions.
Let $\ell$ be a slowly varying function. 
Then $\overline F$ belongs to the de Haan class $\Pi_\ell$
with index $-c$ if
\begin{equation} \label{eq:deHaan}
\lim_{x \to \infty} 
\frac{\overline F(\lambda x) - \overline F(x)}{\ell (x)}
= -c \log \lambda.
\end{equation}
The de Haan class $\Pi$ consists of all functions $f$ such 
that $f \in \Pi_\ell$ with nonzero index for some slowly 
varying $\ell$. We note that $\Pi$ is a proper subclass of the 
slowly varying functions.
For properties on the de Haan class we refer to
Chapter 3 in \cite{BGT} and Appendix B in de Haan and 
Ferreira \cite{deHaan}.

The main novelty in the next result is the relation 
between the two truncations, $h_\beta$ and $V_\beta$.

\begin{theorem}
Assume \eqref{eq:beta-ass} for some $\beta > 0$. 
For any $\rho \in (0,\beta)$ the following are equivalent:
\begin{align}
& h_\beta \in \mathcal{RV}_\rho; \label{eq:h-rv}\\ 
& V_\beta \in \mathcal{RV}_\rho; \label{eq:V-rv} \\
& \overline F \in \mathcal{RV}_{\rho- \beta}; \label{eq:F-rv} \\
& \lim_{x \to \infty} 
\frac{\overline F(x) x^\beta}{h_\beta(x)} = \frac{\rho}{\beta}; 
\label{eq:lim1} \\
& \lim_{x \to \infty} \frac{V_{\beta}(x)}{h_\beta(x)} = 
1 - \frac{\rho}{\beta}. \label{eq:lim2}
\end{align}

If $\rho = 0$ then conditions \eqref{eq:h-rv}, \eqref{eq:V-rv},
\eqref{eq:lim1}, and \eqref{eq:lim2} are equivalent,
and \eqref{eq:F-rv} implies each of them.

If $\rho = \beta$ then conditions \eqref{eq:h-rv}, \eqref{eq:F-rv},
\eqref{eq:lim1}, and \eqref{eq:lim2} are equivalent, and 
\eqref{eq:V-rv} implies each of them. Furthermore,
\eqref{eq:V-rv} is equivalent to $\overline F \in \Pi$.
\end{theorem}
%
%
%

\begin{remark}
For $\beta = 1$ and $\rho = 0$
the result follows from Rogozin's theorem \cite[Theorem 8.8.1]{BGT}, and
\eqref{eq:V-rv} is further equivalent to the relative stability of the
corresponding random walk, see \cite[Section 8.8]{BGT}. 

Another important
special case is $\beta = 2$, when $V_2$ is the truncated second moment.
A distribution belongs to the domain of attraction of the normal law
if and only if $V_2$ is slowly varying, see \cite[Section 8.3]{BGT}.
In the corresponding small time setup (asymptotics at 0), 
for $\beta = 2$, $\rho = 0$, among other results Theorem 2.4 by 
Maller and Mason
\cite{MallerMason10} states that both \eqref{eq:h-rv} and
\eqref{eq:V-rv} follows from \eqref{eq:lim1}.
\end{remark}

\begin{remark}
For $\rho = 0$ in Theorem 3 \cite{Seneta19} Seneta showed that
\eqref{eq:V-rv}, \eqref{eq:lim1} and \eqref{eq:lim2} are 
equivalent, and any of them implies \eqref{eq:h-rv}.
The implication \eqref{eq:h-rv} $\Rightarrow$ \eqref{eq:V-rv}
is the conjecture by Seneta \cite[(26) on p.1653]{Seneta19}.
Moreover, in \cite[Theorem 3]{Seneta19} it is also shown that 
$h_\beta$ is in fact \emph{slowly varying in the Zygmund sense},
which is, by a result of Bojani\'c and Karamata \cite[Theorem 1.5.5]{BGT},
equivalent to $h_\beta$ being normalized slowly varying.
\end{remark}

\begin{remark}
For $\rho \in (0,\beta)$ the functions $x^\beta \overline F(x)$, 
$h_{\beta}(x)$, and $V_\beta(x)$ are all asymptotically equivalent
up to a strictly positive finite constant factor. The borderline 
cases $\rho =0$ and $\rho = \beta$ are somewhat different.

For $\rho = 0$ we have $h_\beta(x) \sim V_\beta(x)$, and 
$x^{\beta} \overline F(x) = o(h_\beta(x))$. Moreover, the 
slow variation of $h_\beta$ does not imply the slow variation 
of $x^\beta \overline F(x)$. Indeed, if 
$u(x):= x^{\beta} \overline F(x)$
is a logarithmically periodic function for large $x$, meaning that 
for some $p > 1$ we have $u(x) = u(px)$, then
$\int_0^x u(y) y^{-1} \dd y$ is slowly varying.
See Lemma 2.3 by Kevei \cite{Kevei}, or in a more general setting
Proposition 6.7 by Buldygin et al.~\cite{BIKS}.
The simplest such example for $\beta = 1$
is the classical St.~Petersburg distribution, with
distribution function 
\[ 
F(x) = 
\begin{cases}
1 - 2^{-\lfloor \log_2 x \rfloor} = 1 - \frac{2^{\{ \log_2 x \}}}{x},
& \text{for } x \geq 2, \\
0, & \text{otherwise},
\end{cases}
\]
where $\log_2$ stands for the logarithm with base 2, $\lfloor \cdot \rfloor$
is the (lower) integer part, and $\{ \cdot \}$ is the fractional part;
see \cite[Section 8.8.2]{BGT}, or \cite[Section VII.7]{Feller}.

On the other hand, if $\rho = \beta$ then
$h_\beta \in \mathcal{RV}_\beta$ implies that $\overline F(x)$
is slowly varying and $h_\beta(x) \sim \overline F(x) x^\beta$.
Then $V_\beta(x) = o(h_\beta(x))$, and
\begin{equation} \label{eq:Vbetarho}
x^{-\beta} V_\beta (x) 
= x^{-\beta} \int_0^x \beta y^{\beta-1} 
\overline F(y)  \dd y - \overline F(x), 
\end{equation}
which is not necessarily slowly varying.
In fact, by Theorem 3.7.1 in \cite{BGT} (a version of de Haan's theorem)
it is slowly varying if and only if $\overline F$
belongs to the de Haan class $\Pi$.
\end{remark}

\section{Proof}

\emph{Equivalence \eqref{eq:lim1} $\Leftrightarrow$ \eqref{eq:lim2}}.
Clearly, \eqref{eq:int-parts} implies that 
\eqref{eq:lim1} and \eqref{eq:lim2} are equivalent for any 
$\rho \in [0,\beta]$.
\medskip

\emph{Implications \eqref{eq:h-rv} $\Rightarrow$ 
\eqref{eq:V-rv}, \eqref{eq:F-rv}, \eqref{eq:lim1}}.
Assume \eqref{eq:h-rv}.
Then for $\lambda > 1$
\[
\begin{split}
h_\beta(\lambda x) - h_\beta(x) & = 
\beta \int_{x}^{\lambda x} y^{\beta-1} \overline F(y) \dd y \\
& \geq \overline F(\lambda x) \left( (\lambda x)^\beta - x^\beta \right) \\
& = \overline F(\lambda x) (\lambda x)^\beta ( 1 - \lambda^{-\beta}).
\end{split}
\]
Dividing both sides by $h_\beta(\lambda x)$, taking limits as $x \to \infty$
and using that $h_\beta$ is regularly varying we obtain
\[
\limsup_{x \to \infty} \frac{\overline F(\lambda x) (\lambda x)^\beta}
{h_\beta(\lambda x)} \leq \frac{1 - \lambda^{-\rho}}{1 - \lambda^{-\beta}}.
\]
As $\lambda \downarrow 1$ we have
\begin{equation} \label{eq:limsup}
\limsup_{x \to \infty} \frac{\overline F(x) x^\beta}{h_\beta(x)}
\leq \frac{\rho}{\beta}.
\end{equation}

Similarly we derive the lower bound. 
For $\lambda > 1$
\[
\begin{split}
h_\beta(\lambda x) - h_\beta(x) 
& \leq \overline F(x) \left( (\lambda x)^\beta - x^\beta \right) \\
& = \overline F(x) x^\beta ( \lambda^\beta - 1).
\end{split}
\]
Dividing both sides by $h_\beta(x)$ and taking limits as $x \to \infty$
we obtain
\[
\liminf_{x \to \infty} \frac{\overline F(x) x^\beta}
{h_\beta(x)} \geq \frac{\lambda^{\rho} - 1}{\lambda^{\beta}-1}.
\]
As $\lambda \downarrow 1$ we have
\begin{equation} \label{eq:liminf}
\liminf_{x \to \infty} \frac{\overline F(x) x^\beta}{h_\beta(x)}
\geq \frac{\rho}{\beta}.
\end{equation}
Combining \eqref{eq:limsup} and \eqref{eq:liminf} 
we obtain \eqref{eq:lim1}. Furthermore, if $\rho < \beta$ then
by \eqref{eq:lim2} $V_\beta(x) \sim h_\beta(x) (1 - \rho/\beta)$, thus
\eqref{eq:V-rv} also follows. While, if $\rho > 0$ then \eqref{eq:lim1}
implies \eqref{eq:F-rv}.
\medskip

\emph{Implications \eqref{eq:V-rv} $\Rightarrow$ \eqref{eq:h-rv},
\eqref{eq:lim1}}.
Assume \eqref{eq:V-rv}.
Since $\int_{[0,\infty)} y^\beta \dd F(y) = \infty$,
we may apply Theorem 8.1.2 in \cite{BGT} with $\alpha = 0$,
(see also Feller \cite[Section VIII.9]{Feller}) 
and we have
\begin{equation} \label{eq:gamma-lim}
\lim_{x \to \infty} \frac{x^\beta \overline F(x)}{V_{\beta}(x)} 
= \gamma \in [0,\infty],
\end{equation}
and there exists $p \in [0,\beta]$ and 
a slowly varying function $\ell \in \mathcal{RV}_0$ such that
\begin{equation} \label{eq:p}
\gamma = \frac{\beta - p}{p}, \quad 
\frac{V_{\beta}(x)}{x^{\beta - p} \ell(x)} \to p, \quad 
\frac{x^p \overline F(x)}{\ell(x)} \to  \beta - p.
\end{equation}

If $\gamma \in (0,\infty)$ in \eqref{eq:gamma-lim} then
the second convergence in \eqref{eq:p} implies
$p = \beta - \rho \in (0,\beta)$, thus 
$\gamma = \rho / ( \beta - \rho)$. 
Therefore, by \eqref{eq:p}
\[
\begin{split}
\frac{h_\beta(x)}{V_\beta(x)}  = 1 + 
\frac{x^\beta \overline F(x)}{V_\beta(x)} 
\to 1 + \frac{\rho}{\beta - \rho} 
= \frac{\beta}{\beta - \rho},
\end{split}
\]
proving both \eqref{eq:h-rv} and \eqref{eq:lim1}.

If $\gamma = 0$ then by \eqref{eq:p} $p = \beta$, 
which implies $\rho = 0$, 
$V_\beta (x) \sim \beta \ell(x)$,
and $x^\beta \overline F(x) = o(\ell(x))$. Thus 
$h_\beta(x) \sim V_\beta(x)$ by \eqref{eq:int-parts}, 
implying \eqref{eq:h-rv} and \eqref{eq:lim1}.

If $\gamma = \infty$ then $p = 0$, which implies $\rho = \beta$, 
$\overline F(x) \sim \beta \ell(x)$, and 
$V_\beta (x) = o(x^{\beta} \ell(x))$. Therefore, $h_\beta(x) 
\sim x^\beta \overline F(x)$, in particular, \eqref{eq:h-rv} 
and \eqref{eq:lim1} holds.
\medskip

\emph{Implication \eqref{eq:F-rv} $\Rightarrow$ \eqref{eq:h-rv}}.
For any $\rho > 0$ this is an immediate consequence of 
Karamata's theorem (\cite[Proposition 1.5.8]{BGT}). If $\rho = 0$
then by Proposition 1.5.9a in \cite{BGT} $h_\beta$ is slowly varying,
and $h_\beta(x) / (x^{\beta} \overline F(x) ) \to \infty$.
\medskip

\emph{Implication \eqref{eq:lim1} $\Rightarrow$ \eqref{eq:h-rv}}.
Assume \eqref{eq:lim1}. Then, since 
$h_\beta'(x) = \beta x^{\beta -1} \overline F(x)$ Lebesgue almost everywhere,
we have
\begin{equation} \label{eq:delta}
\delta(x) := \frac{x h_\beta'(x)}{h_\beta(x)} \to \rho.
\end{equation}
Therefore for some $A > 0$, $B> 0$
\[
h_\beta(x) = A \exp \left( \int_B^x \frac{\delta(y)}{y} \dd y \right)
= A  B^{-\rho} \, x^{\rho} 
\exp \left( \int_B^x \frac{\delta(y)-\rho}{y} \dd y \right).
\]
By the representation theorem and \eqref{eq:delta}
the second factor is a slowly varying function (in fact 
it is normalized slowly varying), showing that 
$h_\beta \in \mathcal{RV}_\rho$, i.e.~\eqref{eq:h-rv} holds.
\medskip

The only remaining part to prove is that for 
$\rho = \beta$ conditions $V_\beta \in \mathcal{RV}_\beta$
and $\overline F \in \Pi$ are equivalent.
By \eqref{eq:Vbetarho} and Theorem 3.7.1 in \cite{BGT} 
(a version of de Haan's theorem) we see that
$x^{-\beta} V_\beta(x) \sim c \ell(x) / \beta$ for some $c > 0$
and slowly varying $\ell$ if and only if \eqref{eq:deHaan} holds.
The proof is complete.

\bigskip

\noindent
\textbf{Acknowledgements.}
This research was partially supported 
by the J\'{a}nos Bolyai Research Scholarship of the 
Hungarian Academy of Sciences, 
by the NKFIH grant FK124141, 
and by the
Ministry of Human Capacities, Hungary grant 
TUDFO/47138-1/2019-ITM.


\end{document}